%% file: Near_Counterexamples_to_Weils_Converse_Theorem.tex
\documentclass[11pt, a4paper]{amsart}
\usepackage{amsmath}%
\usepackage{amsfonts}%
\usepackage{amssymb}%
\usepackage{graphicx}
\usepackage{amsthm}
\usepackage{amssymb}
\usepackage{tikz}
\usepackage{tikz-cd}
\usetikzlibrary{matrix,arrows}
\usepackage[english, ngerman]{babel}
\usepackage[utf8]{inputenc}
\usepackage{hyperref}
\usepackage{cite}
\usepackage{footnote}
\usepackage{stmaryrd}
\usepackage[a4paper,margin=3cm]{geometry}
\usepackage{verbatim}
\usepackage{enumerate}

\newtheorem{thm}{Theorem}
\theoremstyle{plain}

\newtheorem{cor}[thm]{Corollary}

\theoremstyle{remark}
\newtheorem{rem}[thm]{Remark}
\newtheorem*{acknowledgements}{Acknowledgements}

\renewcommand{\phi}{\varphi}

\newcommand{\BC}{{\mathbb{C}}}

\newcommand{\BH}{{\mathbb{H}}}

\newcommand{\BN}{{\mathbb{N}}}

\newcommand{\BR}{{\mathbb{R}}}

\newcommand{\BZ}{{\mathbb{Z}}}

\newcommand{\CQ}{{\mathcal Q}}

\newcommand{\nin}{\notin}

\renewcommand{\mod}{\mathop{\rm mod}\nolimits}

\newcommand{\SL}[1]{\mathop{\rm SL}_{#1} \nolimits}

\renewcommand{\Re}{\mathop{\rm Re}\nolimits}

\newcommand{\quotient}[2]{
        \mathchoice
            {
                \text{\raise1ex\hbox{$#1$}\Big/\lower1ex\hbox{$#2$}}%
            }
            {
                #1\,/\,#2
            }
            {
                #1\,/\,#2
            }
            {
                #1\,/\,#2
            }
    }
		
\newcommand{\rquotient}[2]{
        \mathchoice
            {
                \text{\lower1ex\hbox{$#1$}\Big \backslash \raise01ex\hbox{$#2$}}%
            }
            {
                #1\,\backslash\,#2
            }
            {
                #1\,\backslash\,#2
            }
            {
                #1\,\backslash\,#2
            }
    }
		
\newcommand{\lrquotient}[3]{
        \mathchoice
            {
                \text{\lower1ex\hbox{$#1$}\Big \backslash \raise01ex\hbox{$#2$}\Big/\lower1ex\hbox{$#3$}}%
            }
            {
                #1\,\backslash\,#2\,/\,#3
            }
            {
                #1\,\backslash\,#2\,/\,#3
            }
            {
                #1\,\backslash\,#2\,/\,#3
            }
    }

\begin{document}
\selectlanguage{english}

\bibliographystyle{plain}

\title{Near Counterexamples to Weil's Converse Theorem}
\author{Raphael S. Steiner}
\address{Department of Mathematics, University of Bristol, Bristol BS8 1TW, UK}%
\email{raphael.steiner@bristol.ac.uk}%




\begin{abstract} We show that in Weil's converse theorem the functional equations of multiplicative twists for at least the first $\sqrt{\frac{p-24}{3}}$ moduli are needed in order to prove the modularity for $\Gamma_0(p)$.
\end{abstract}
\maketitle

\input{Introduction}

\begin{acknowledgements} I would like to thank Andrew Booker for asking me the question whether the modularity of $|f|$ implies the modularity of $f$, which eventually led to the above observation, and for interesting discussions as well as useful comments on earlier versions of this paper.
\end{acknowledgements}
\input{Notation}

\input{Proof}

\bibliography{Bibliography}
\end{document}

%% file: Introduction.tex
\section{Introduction}

Let $\SL{2}(\BR)$ act on the upper half-plane $\BH$ by M\"obius-transformations. That is, we define as usual
$$
\gamma z = \frac{az+b}{cz+d} \text{ for } \gamma = \begin{pmatrix} a & b \\ c & d \end{pmatrix} \in \SL{2}(\BR).
$$
For such a matrix we further define the cocycle $j(\gamma,z)=cz+d$ and an action on the holomorphic functions $f: \BH \to \BC$ on the upper-half-plane as $(f|_k\gamma)(z)=j(\gamma,z)^{-k}f(\gamma z)$. For our convenience let us denote the following special matrices
$$
S= \begin{pmatrix} 1 & 1 \\ 0 & 1 \end{pmatrix}, \quad \quad T= \begin{pmatrix} 0 & -1 \\ 1 & 0 \end{pmatrix}, \quad \quad W_N=\begin{pmatrix} 0 & -N^{-\frac{1}{2}} \\ N^{\frac{1}{2}} & 0  \end{pmatrix}.
$$
Let $(a_m)_{m \in \BN},(b_m)_{m \in \BN}$ be two complex sequences that satisfy $|a_m|,|b_m|=O(m^{\sigma})$ for some $\sigma >0 $, where we follow Landau's big $O$ convention. Associated to these sequences we define the two holomorphic functions on the upper half-plane
$$
f(z)= \sum_{m \ge 1} a_m e(mz) \text{ and } g(z)= \sum_{m \ge 1}b_m e(mz),
$$
where $e(z)=e^{2 \pi i z}$. Let us further define their associated Dirichlet series and completions thereof
$$\begin{aligned}
L(f,s)&=\sum_{m \ge 1}a_m m^{-s}, && L(g,s)= \sum_{m \ge 1}b_m m^{-s}, \\
\Lambda(f,s)&=(2 \pi)^{-s} \Gamma(s) L(f,s), && \Lambda(g,s)=(2 \pi)^{-s} \Gamma(s) L(g,s),
\end{aligned}$$
where $\Gamma(s)$ is the Euler Gamma function. These are \emph{a priori} only defined for $\Re(s) > \sigma +1$. For $k\in \BZ$, Hecke \cite{HeckeConverse} proved the equivalence of the following two statements.
\begin{enumerate}[(i)]
	\item $f,g $ are cusp forms of weight $k$ for $\SL{2}(\BZ)$, and $(f|_k T)(z)=g(z)$.
	\item The functions $\Lambda(f,s)$ and $\Lambda(g,s)$ admit a holomorphic continuation to the whole complex plane, are bounded in any vertical strip and satisfy the functional equation
	$$
	\Lambda(f,s)=i^k\Lambda(g,k-s).
	$$
\end{enumerate}
An equivalent statement for
$$
\Gamma_0(N) = \left \{ \begin{pmatrix} a & b \\ c & d \end{pmatrix} \in \SL{2}(\BZ)  \bigg| \, N|c \right \}
$$
is much harder. Especially the implication going from the functional equations to proving that $f$ and $g$ must be cusp forms requires more information. Weil \cite{WeilConverse} resolved this issue by assuming further functional equations coming from multiplicative twists. For this define for a primitive character $\psi$ modulo $q$ with $(q,N)=1$ the Dirichlet series and their respective completions
$$\begin{aligned}
L(f, \psi ,s)&=\sum_{m \ge 1}a_m \psi(m) m^{-s}, && L(g,\psi,s)= \sum_{m \ge 1}b_m \psi(m) m^{-s}, \\
\Lambda(f, \psi,s)&=(2 \pi)^{-s} \Gamma(s) L(f,\psi,s), && \Lambda(g,\psi,s)=(2 \pi)^{-s} \Gamma(s) L(g,\psi,s).
\end{aligned}$$
Assume the holomorphic continuation of $\Lambda(f,\psi,s),\Lambda(g,\psi,s)$ to the whole complex plane and boundedness in any vertical strip for every primitive character $\psi$ modulo $q$ with $(q,MN)=1$, where $M$ is some fixed integer. Furthermore assume all the functional equations
$$
\Lambda(f,\psi,s)= i^k \chi(q) \psi(N) \frac{\tau(\psi)^2}{q} (Nq^2)^{\frac{k}{2}-s} \Lambda(g,\overline{\psi},k-s),
$$
for some fixed Dirchlet character $\chi \mod(N)$, where $ \tau(\psi)$ is the Gauss sum of $\psi$. In this case Weil \cite{WeilConverse} was able to prove that $f$ respectively $g$ are cusp forms of weight $k$ for $\Gamma_0(N)$ with character $\chi$ respectively $\overline{\chi}$ and $(f|_k W_N)(z)=g(z)$. Khoai \cite{Khoai} later refined the number of functional equations of twists needed. He proved that the primitive twists $\psi$ of modulus $q<N^2$ with $(q,N)=1$ suffice to come to the same conclusion. Moreover if $N=p^r$ is a power of a prime, then already $q < p^r$ suffices. In this paper we are able to show that for $N=p$ a prime, then one needs at least the multiplicative twists for the first $(\frac{1}{\sqrt{3}}+o(1))\sqrt{p}$ moduli.

\begin{thm} Let $p \ge 29$ be a prime, $\chi$ an even Dirichlet character modulo $p$ and $k \ge 16$ an even integer. There exists two complex sequences $(a_m)_{m \in \BN},(b_m)_{m \in \BN}$ with $|a_m|,|b_m|=O(m^{\sigma})$ for some $\sigma > 0$ which satisfies the following two properties:
\begin{enumerate}[(i)]
	\item Both $f(z)=\sum_{m \ge 1 }a_m e(mz)$ and $g(z)=\sum_{m \ge 1}b_m e(mz) \nin \bigcup_{l,N} S_l(\Gamma(N))$, that is to say $f$ and $g$ are not classical cusp forms of any weight for any congruence subgroup,
	\item for any primitive character $\psi$ modulo $q$ with $(q,p)=1$ the functions $\Lambda(f,\psi,s),\Lambda(g,\psi,s)$ can be holomorphically continued to the whole complex plane, are bounded in any vertical strip and satisfying the functional equations
$$
\Lambda(f,\psi,s)= i^k\chi(q)\psi(p)\frac{\tau(\psi)^2}{q} (pq^2)^{\frac{k}{2}-s} \Lambda(g,\overline{\psi},k-s), \quad \forall q \le \sqrt{\frac{p-24}{3}}.
$$ 
\end{enumerate}
\label{thm:mult}
\end{thm}
In the context of converse theorems for classical modular forms it is more natural to look at additive twists rather than multiplicative ones. For this purpose we define for $(a,q)=1$:

$$\begin{aligned}
L\left(f, \frac{a}{q} ,s\right)&=\sum_{m \ge 1}a_m e\left(\frac{am}{q}\right) m^{-s}, && L\left(g,\frac{a}{q},s\right)= \sum_{m \ge 1}b_m e\left(\frac{am}{q}\right) m^{-s}, \\
\Lambda\left(f, \frac{a}{q},s\right)&=(2 \pi)^{-s} \Gamma(s) L\left(f,\frac{a}{q},s\right), && \Lambda\left(g,\frac{a}{q},s\right)=(2 \pi)^{-s} \Gamma(s) L\left(g,\frac{a}{q},s\right).
\end{aligned}$$
For additive twists we are able to give the following converse theorem for prime level $p$.

\begin{thm} Let $p>3$ be a prime. There exist a (computable) set $\CQ$ of $2\lfloor \frac{p}{12} \rfloor+3$ natural numbers $1 \le q \le p-2$ for which the following holds. Given any two complex sequences $(a_m)_{m\in \BN}$ and $(b_m)_{m \in \BN}$ with $|a_m|,|b_m|=O(m^{\sigma})$ for some $\sigma > 0$, then the following two statements are equivalent:
\begin{enumerate}[(i)]
	\item $f(z)=\sum_{m \ge 1} a_m e(mz)$ is a cusp form of weight $k$ for $\Gamma_0(p)$ with character $\chi$ and $(f|_k W_p)(z)=g(z)=\sum_{m \ge 1} b_m e(mz)$,
	\item The functions $\Lambda(f,\frac{-1}{q},s)$ and $\Lambda(g,\frac{(qq_{\star}+1)/p}{q},s)$, where $1\le q_{\star} \le p$ and satisfies $qq_{\star}\equiv -1 \mod(p)$, for every $q	\in \CQ$ can be holomorphically continued to the whole complex plane and are bounded in every vertical strip. Moreover they satisfy the functional equations
	$$\Lambda\left(f, \frac{-1}{q},s\right)=i^k \chi(q) (pq^2)^{\frac{k}{2}-s} \Lambda\left(g,\frac{(qq_{\star}+1)/p}{q},k-s\right), \quad \forall q \in \CQ.$$
\end{enumerate}
\label{thm:add}
\end{thm}

\begin{rem} The number of functional equations of additive twists to go back from \emph{(ii)} to \emph{(i)} in Theorem \ref{thm:add} is essentially optimal as one can construct counterexamples similar to the proof of Theorem \ref{thm:mult} if one assumes at most $2\lfloor \frac{p}{12} \rfloor -2$ additive twists.
\label{rem:thm2}
\end{rem}


Note that this further strengthens Khoai's result, since by using Gauss sums one may reduce down to about $p/6$ moduli for which one needs to consider multiplicative twists.

The proof of Theorem \ref{thm:add} is straightforward and relies on Hecke's converse theorem \cite{HeckeConverse} and on a result of Rademacher \cite{Rademacher29} on the generators of $\Gamma_0(p)$. The proof of Theorem \ref{thm:mult} however relies on the interesting observation that the modularity of $|f|$ for a finite index subgroup $\Gamma \subseteq \SL{2}(\BZ)$ does not always imply the modularity of $f$ on some congruence subgroup $\Gamma(N)$, despite this being the case for the full modular group $\Gamma=\SL{2}(\BZ)$ where it would follow that $f$ is modular on $\Gamma(12)$. In fact for $\Gamma=\Gamma_0(p)$ with $p \to \infty$ there is a vast amount of freedom which we shall exploit.

%% file: Notation.tex
\section{Notation and Preliminaries} Let $\Gamma$ be a finite index subgroup of $\SL{2}(\BZ)$ that contains $-I$. We call a character $\upsilon : \Gamma \to S^1$ that satisfies $\upsilon(-I)=(-1)^k$ a multiplier system for $\Gamma$. We call a holomorphic function $f$ on the upper half plane modular of weight $k$ for $\Gamma$ with respect to $\upsilon$ if it satisfies $(f|_k \gamma)(z)=\upsilon(\gamma) f(z)$ for every $\gamma \in \Gamma$. Such a function $f$ has an expansion of the shape
$$
(f|_k \tau)(z)=\sum_{m=-\infty}^{\infty} a_m e\left( \frac{m+\kappa_{\tau}}{n_{\tau}} z \right)
$$
for every $\tau \in \SL{2}(\BZ)$, where $n_{\tau}$ is the width of the cusp $\tau \infty$ and $\kappa_{\tau}$ is the cusp parameter at the cusp $\tau \infty$. They are both independent of the choice of representative of $\tau \infty \ \mod(\Gamma)$. The former is characterised by being the smallest natural number $n$ such that $\tau S^n \tau^{-1} \in \Gamma$ and the latter is characterised by $e(\kappa_{\tau})=\upsilon(\tau S^{n_\tau} \tau^{-1})$ and $\kappa_{\tau} \in [0,1)$. We say $f$ is a modular form of weight $k$ for $\Gamma$ with respect to $\upsilon$ if we can restrict the summation to $m+\kappa_{\tau} \ge 0$ for every $\tau \in \SL{2}(\BZ)$ and moreover we say $f$ is a cusp form of weight $k$ for $\Gamma$ with respect to $\upsilon$ if the summation can be restricted to $m+\kappa_{\tau}>0$ for every $\tau \in \SL{2}(\BZ)$. For a more detailed treatment on modular forms (of arbitrary real weight) with respect to an arbitrary multiplier system we refer the reader to \cite{MFaF}.

From now on let $k$ denote an even integer and let $\Gamma= \Gamma_0(p)$, where $p>3$ is a prime. The multiplier systems for $\Gamma_0(p)$ include all even Dirichlet characters $\chi: (\BZ \slash p \BZ)^{\times} \to S^1$; they are given by
$$
\upsilon_{\chi}(\gamma) = \chi(d) \text{ if } \gamma=\begin{pmatrix} a & b \\ c & d \end{pmatrix} \in \Gamma_0(p),
$$
but in fact there are many more multiplier systems. They are in one-to-one correspondence with group homomorphisms
$$
\upsilon: \quotient{\Gamma_{0}(p)}{\pm[\Gamma_{0}(p),\Gamma_{0}(p)]} \to S^1.
$$
We can characterize this quotient completely by a theorem of Rademacher \cite{Rademacher29}.

\begin{thm}[Rademacher] Let $p>3$ be a prime, then we have
$$
\quotient{\Gamma_{0}(p)}{\pm I} \cong F_{l-2a-2b} \ast  \left(\BZ \slash 2 \BZ \ast \BZ \slash 2 \BZ \right)^{a} \ast \left(\BZ \slash 3 \BZ \ast \BZ \slash 3 \BZ \right)^{b},
$$
where $l=2 \lfloor \frac{p}{12} \rfloor + 3$, $a=0$ unless $p \equiv 1 \mod(4)$ in which case $a=1$ and $b=0$ unless $p \equiv 1 \mod(3)$ in which case $b=1$ and $F_{l-2a-2b}$ is the free group with $l-2a-2b$ generators. The isomorphism comes from a set of free generators of $\Gamma_0(p) \slash \{\pm I\}$, which are given by $S$ and some matrices of the shape
$$
V_q= \begin{pmatrix} -q_{\star} & -1 \\ qq_{\star}+1 & q \end{pmatrix},
$$
with $2 \le q \le p-2$.
\label{thm:rade}
\end{thm}
As a corollary we find
\begin{cor} Let $p>3$ be a prime, then we have
$$
\quotient{\Gamma_{0}(p)}{\pm[\Gamma_{0}(p),\Gamma_{0}(p)]} \cong \BZ^{l-2a-2b} \times \left(\BZ \slash 2 \BZ \right)^{2a} \times \left(\BZ \slash 3 \BZ \right)^{2b},
$$
where $l=2 \lfloor \frac{p}{12} \rfloor + 3$, $a=0$ unless $p \equiv 1 \mod(4)$ in which case $a=1$ and $b=0$ unless $p \equiv 1 \mod(3)$ in which case $b=1$. 
\label{cor:rade}
\end{cor}
Thus when $p$ becomes large we get plenty of freedom, which we may use to satisfy certain equations.

%% file: Proof.tex
\section{Proof of the Theorems}

Let $f(z)=\sum_{m > 0} a_m e(mz)$ be a cusp form of even weight $k>0$ with respect to a multiplier system $\upsilon$ on $\Gamma_0(p)$ with cusp parameters $\kappa_{I}=\kappa_T=0$ and let us denote
$$
g(z)=(f|_k W_p)(z)=\sum_{m > 0 }b_m e(mz),
$$
which is another cusp form with respect to a conjugated multiplier system of $\upsilon$.

\begin{rem} The assumption on $\kappa_T$ is necessary as the example $p=13$ shows, since there one has $TS^{13}T^{-1}=V_{10}^{-2}V_8^{-1}V_5^{-1}V_4^{-2}S^{-1}$. Although it is interesting to note that given $\kappa_I=0$ then $\upsilon(TS^pT^{-1})$ is always a $6$-th root of unity and if $p \equiv 11 \mod(12)$ then it is equal to $1$. This can be seen from the relation $TS^pT^{-1}=V_{p-1}S^{-1}$ and going through Rademacher's elimination process \cite{Rademacher29}.
\end{rem}

Let $q \in \BN$ with $q\neq p$ and $0 \le a < q$ with $(a,q)=1$. Let $B,D$ be two integers that satisfy the relation $qD+apB=1$, which has solutions by B\'ezout. We have the matrix identity
$$
\begin{pmatrix} 1 & \frac{a}{q} \\ 0 & 1 \end{pmatrix} = \begin{pmatrix} D & a \\ -pB & q \end{pmatrix} W_p \begin{pmatrix} p^{\frac{1}{2}}B & p^{-\frac{1}{2}}q^{-1} \\ -p^{\frac{1}{2}}q & 0 \end{pmatrix}
$$
and thus we get the identity
\begin{equation}
\left(f \bigg|_k \begin{pmatrix} 1 & \frac{a}{q} \\ 0 & 1 \end{pmatrix} \right)(z) = \upsilon \left( \begin{pmatrix} D & a \\ -pB & q \end{pmatrix} \right) \left( g \bigg|_k \begin{pmatrix} p^{\frac{1}{2}}B & p^{-\frac{1}{2}}q^{-1} \\ -p^{\frac{1}{2}}q & 0 \end{pmatrix} \right)(z).
\label{eq:modrelation}
\end{equation}
By definition this is just
\begin{equation}\begin{aligned}
\sum_{m > 0} a_m & e\left( \frac{am}{q} \right) e(mz) \\
&= (-1)^k \upsilon \left( \begin{pmatrix} D & a \\ -pB & q \end{pmatrix} \right) p^{-\frac{k}{2}} q^{-k} z^{-k} \sum_{m > 0} b_m e\left( -\frac{Bm}{q} \right) e\left( - \frac{m}{pq^2z} \right).
\end{aligned}
\label{eq:modularity}
\end{equation}
Now $f|_k \bigl(\begin{smallmatrix} 1 & a/q \\ 0 & 1 \end{smallmatrix}\bigr)$ and $g|_k \bigl(\begin{smallmatrix} 1 & -B/q \\ 0 & 1 \end{smallmatrix}\bigr)$ are modular forms on $\Gamma(pq^2)$ for some multiplier system, respectively, and thus by Hecke's original proof $\Lambda(f,\frac{a}{q},s)$ and $\Lambda(g,\frac{-B}{q},s)$ have a holomorphic continuation to the whole complex plane and are bounded in every vertical strip. Now by Bochner \cite{Boch51} the modular relation \eqref{eq:modularity} is equivalent to the functional equation
\begin{equation}
\Lambda\left(f,\frac{a}{q},s\right) = i^k \upsilon \left( \begin{pmatrix} D & a \\ -pB & q \end{pmatrix} \right) (pq^2)^{\frac{k}{2}-s} \Lambda\left(g,\frac{-B}{q},k-s\right).
\label{eq:fqeq}
\end{equation}
This shows one direction in Theorem \ref{thm:add}. For the other direction we need to fix our set $\CQ$. We take it to consist of $1$ and those $q$ for which $V_q$ is needed to generate $\quotient{\Gamma_0(p)}{\pm I}$ in Theorem \ref{thm:rade}. We make use of the equivalence of \eqref{eq:fqeq} and \eqref{eq:modrelation}, where $\upsilon \bigl( \bigl(  \begin{smallmatrix}D&a\\ -pB & q \end{smallmatrix} \bigr)\bigr)$ is to be regarded as any fixed constant of modulus 1. We first set $(a,q,B,D)=(-1,1,-1,1-p)$ and make use of the functional equation for $1\in \CQ$. This yields
$$
\left(f \bigg|_k \begin{pmatrix} 1 & -1 \\ 0 & 1 \end{pmatrix} \right)(z) = \chi(1) \left( g \bigg|_k \begin{pmatrix} -1 & -1 \\ 0 & -1 \end{pmatrix} W_p \right)(z)
$$
or equivalently $f(z)=(g|_kW_p)(z)$ and consequently $(f|_kW_p)(z)=g(z)$ after applying $|_kW_p$ to both sides.
We then further find that \eqref{eq:modrelation} is equivalent to
$$
\left(f \bigg|_k \begin{pmatrix} D & a \\ -pB & q \end{pmatrix} \right)(z) = \upsilon \left( \begin{pmatrix} D & a \\ -pB & q \end{pmatrix} \right) \left( g |_k W_p \right)(z).
$$
Using this with $(a,q,B,D)=(1,q,-(qq_{\star}+1)/p,-q_{\star})$, where $q \in \CQ$, and the functional equation for $q$ we find
$$
\left(f \bigg|_k \begin{pmatrix} -q_{\star} & 1 \\ qq_{\star}+1 & q \end{pmatrix} \right)(z) = \chi(q) f(z).
$$
Combining this with Theorem \ref{thm:rade} and the trivial fact that $(f|_k S)(z)=f(z)$ we find that $f$ is modular of weight $k$ for $\Gamma_0(p)$ with character $\chi$. The fact that $f$ is indeed a cusp form comes from the expansions at the cusps $\infty$ and $0$ which are given by the definition of $f$ respectively $g$.\\

Now onto the proof of Theorem \ref{thm:mult}. Let $\upsilon,f$ and $g$ be defined as in the beginning of this section. The restriction that $\kappa_I=\kappa_T=0$ implies that the dependence on $B$ in both
$$
\upsilon\left( \begin{pmatrix} D & a \\ -pB & q \end{pmatrix} \right) \text{ and } \Lambda\left( g, \frac{-B}{q} , s \right)
$$
is only on $B \mod(q)$ (and thus only depends on $q$ and $a \mod q$).
We further make the assumption that
$$
\upsilon\left( \begin{pmatrix} D & a \\ -pB & q \end{pmatrix} \right) = \chi(q)
$$
for $(a,q)=1$ and $1\le q \le Q$, where $\chi$ is a Dirichlet character modulo $p$. That is to say we want $\upsilon$ to pretend to be the Dirichlet character $\chi$ for small values of $q$. 

Let now $\psi$ be a primitive Dirichlet character modulo $q$ with $(q,p)=1$ and $1\le q \le Q$. Then we have
$$\begin{aligned}
\Lambda(f,\psi,s) &= \frac{1}{\tau(\overline{\psi})}  \sideset{}{'} \sum_{a \mod(q)}\overline{\psi}(a)\Lambda\left(f,\frac{a}{q},s \right) \\
&= \frac{i^k\chi(q)}{\tau(\overline{\psi})} (pq^2)^{\frac{k}{2}-s} \sideset{}{'}\sum_{a \mod(q)} \overline{\psi}(a)  \Lambda\left(g, \frac{-\overline{ap}}{q},k-s \right) \\
&= \frac{i^k\chi(q)\psi(-p)}{\tau(\overline{\psi})} (pq^2)^{\frac{k}{2}-s} \sideset{}{'}\sum_{a \mod(q)} \psi(-\overline{ap})  \Lambda\left(g, \frac{-\overline{ap}}{q},k-s \right) \\
&= i^k\chi(q)\psi(p)\frac{\tau(\psi)^2}{q} (pq^2)^{\frac{k}{2}-s} \Lambda(g,\overline{\psi},k-s),
\end{aligned}$$
where the prime in the summation denotes that we are only summing over $(a,q)=1$ and $\overline{a}$ denotes the multiplicative inverse of $a \mod(q)$. This is the functional equation we would expect for classical modular forms with multiplier system $\upsilon_{\chi}$. In order to prove Theorem \ref{thm:mult} we construct such a multiplier system $\upsilon$ of infinite order as such a multiplier system cannot be trivial on any congruence subgroup. We furthermore have to show that there is a non-zero cusp form with respect to that multiplier system for some weight $k$.

Our assumptions on our multiplier system form a system of linear equations with at most $2+Q^2/2$ equations in $\log \upsilon(\gamma)$, where $\gamma$ runs over our generators in Corollary \ref{cor:rade}, moreover this system of equations admits a solution, namely $\upsilon_{\chi}$. Thus if $2+Q^2/2+5 \le 2 \lfloor \frac{p}{12} \rfloor +3$ or $Q\le \sqrt{(p-24)/3}$ we find that the kernel has dimension at least $5$. Thus we can find an element $\upsilon'$ of infinite order in the kernel which satisfies $\log \upsilon'(\gamma)=0$ for any generator $\gamma$ of finite order and hence $\upsilon = \upsilon' \upsilon_{\chi}$ is a multiplier system which satisfies our requirements. The construction of a non-zero cusp form is now straightforward. As in \cite{MFaF} one can construct Eisenstein series for $k>2$, which are non-zero, say for example $G_I(\Gamma_0(p),4,\upsilon  j^{k};z,0)$, and multiply it with the cusp form $\Delta(z)$ of weight $12$ for $\SL{2}(\BZ)$ to get a cusp form of weight $k+12$ with respect to $\upsilon $ on $\Gamma_0(p)$, which concludes the proof of Theorem \ref{thm:mult}.